

\input amstex.tex
\documentstyle{amsppt}


\input BoxedEPS.tex

\SetepsfEPSFSpecial
\HideDisplacementBoxes
\def\figure#1{\medskip \centerline{\BoxedEPSF{fig#1.eps}}
  \smallskip\centerline{Figure #1}\medskip}

\def\id{\text{\rm id}}
\def\sc{\text{\rm Sc}}

\topmatter
\title Dual decompositions of 4-manifolds II: linear invariants 
\endtitle
\rightheadtext{Dual decompositions II}
\author Frank Quinn\endauthor
\thanks Partially supported by the National Science Foundation\endthanks
\subjclass 57R65 57M25\endsubjclass
\date Revised March 2004\enddate
\address Mathematics, Virginia Tech, Blacksburg VA 24061-0123\endaddress
\email quinn\@math.vt.edu\endemail
\abstract This paper continues the study of decompositions of a smooth 4-manifold into
two handlebodies with handles of index $\leq2$. Part I gave existence results in terms of
spines and chain complexes over the fundamental group of the ambient manifold. Here we
assume that one side of a decomposition has
larger fundamental group, and use this to define algebraic-topological invariants.  These reveal a basic asymmetry in these decompositions: subtle changes on one
side can force algebraic-topologically detectable changes on the other. A solvable iteration
of the basic invariant gives an ``obstruction theory'' using lower commutator quotients. By thinking of a 2-handlebody as essentially determined by the links used as attaching maps for its 2-handles this theory can be thought of a giving ``ambient'' link invariants. The moves used are 
related to the grope cobordism of links developed by  Conant-Teichner, and the Cochran-Orr-Teichner filtration of the link concordance groups. The invariants give algebraically sophisticated  ``finite type'' invariants in the sense of Vassilaev.
\endabstract

\endtopmatter

\head 1. Introduction\endhead
A {\it dual decomposition\/} of a 4-manifold is a description as
a union of two handlebodies, each with handles of index
$\leq 2$. In Part I we gave existence results for these in terms of spines and chain
complexes over the fundamental group of the ambient manifold; here we incorporate refinements of the
fundamental group. Suppose $\Gamma\to \pi_1N$ is a homomorphism, then a ``$\Gamma$ decomposition of $N$''
is a dual decomposition
$N=M\cup W$ together with a factorization
$\pi_1W\to
\Gamma\to\pi_1N$. We study these using linear invariants, particularly the cellular chains of the
$\Gamma$ cover of $W$. We will not consider bilinear invariants such as signatures, that incorporate
intersection information.

There is a basic asymmetry in the behavior of dual decompositions: subtle changes of one side can force
algebraic-topologically detectable changes on the other. We use this to probe the subtle changes. If we
think of $M$ as being primarily determined by the link giving the attaching map of the 2-handles, then we
see invariants of $W$ as giving information about this link. Or if we think of $M$ as a surgery
presentation of the 3-manifold $\partial M$, then invariants of $W$ give information about this 3-manifold.
This information is closely related to classical link and 3-manifold invariants, and to
some  ``quantum'' invariants. 
In terms of links these invariants are local and relative, measuring perturbations of a given link, inside a fixed ambient 4-manifold.
However they are better defined and behaved than most global link invariants. They also have direct
connections to a more mainstream topic than classical link theory, namely handlebody theory of smooth
4-manifolds. 

This paper develops the conceptual and formal structure of linear invariants of $\Gamma$
decompositions. The main geometric results, for instance, show that when properly defined the invariants are all
realized. The theory is somewhat elaborate, so in this paper we are content with a careful development that should clarify the possibilities and limitations. Work on
 computation, explicit examples, or exact relations with other invariants should perhaps wait on a clearer idea of a job to be done. 

We describe the invariants in more detail.  If a
decomposition $N=M\cup W$ is deformed in a way that does not change the $Z[\pi_1N]$ chains of $M$ (up to
chain homotopy), then duality shows that the $Z[\pi_1N]$ chains of $W$ do not change either. On this level
there are no invariants. However if there is a factorization $\pi_1W\to \Gamma\to\pi_1N$ then the
$Z[\Gamma]$ chains may vary  within the constraint that the
$Z[\pi_1N]$ reduction is fixed. The general setting is a
homomorphism of groups
$\Lambda\to \Gamma$, and free based chain complexes over $Z[\Lambda]$ with a fixed
$Z[\Gamma]$ reduction. If $M\cup W$ is a decomposition with $\pi_1W\to\Lambda\to\Gamma\to
\pi_1N$ given then the $\Lambda\to \Gamma$ ``chain invariant''  is 
 defined simply as an appropriate equivalence class of such chain data. By design it is  unchanged by $\Lambda$ deformations.

We explain how chain invariants can be seen as obstructions. Suppose $M_0\cup W_0$ is a decomposition with a
factorization $\pi_1W_0\to \Lambda\to \Gamma\to \pi_1N$. The $\Gamma$-deformation move described in \S3 gives a new
$\Gamma$-decomposition
$M_1\cup W_1$ together with a canonical simple chain equivalence of the chains of $W_0$ and $W_1$, with
$Z[\Gamma]$ coefficients.
 Suppose it should happen that
the canonical factorization $\pi_1W_1\to\Gamma\to \pi_1N$ lifts to $\Lambda$. Then the
$\Lambda$ chains of $W_i$ together with the canonical equivalence over $\Gamma$ provide
the data needed to define the $\Lambda\to \Gamma$ chain invariant. Note that if the
deformation were actually a $\Lambda$ deformation then the complexes would be equivalent
over
$Z[\Lambda]$ and the invariant of $W_1$ would be the same as that of $W_0$. Therefore the
invariant gives an obstruction to
$M_1$ together with this particular lift  $\pi_1W_1\to \Lambda$ to come from a $\Lambda$
deformation. 

The factorization $\pi_1W\to \Lambda$ in the previous paragraph presents a problem. The
chain invariant is not defined unless a factorization exists, and then may depend on which
one is used. However we observe that given
$\pi_1W\to \Gamma$ there is a particular $\Lambda$ for which canonical factorizations
exist. This is $\Gamma_1(W)=\pi_1W/\text{ker}(\pi_1N\to\Gamma)^{(1)}$, where the
superscript ``${(1)}$'' in the denominator denotes the commutator subgroup. $\Gamma_1(W)$
is an abelian extension of $\Gamma$, and the chain invariant in this case
is an analog of the classical abelian link invariants.

The abelian extension construction can be iterated  to get an ``obstruction theory'' for invariants  using
solvable extensions of $\Gamma$. Specifically we define
$\Gamma_n(W)$ to be $\pi_1W$ divided by the $n$-fold iterated commutator of the kernel of
the homomorphism to $\Gamma$. Suppose two decompositions have the same
$n$-fold solvable extension $\Gamma_n(W_0)\simeq \Gamma_n(W_1)$, and the same chains
over this extension. Then there is a canonical isomorphism $\Gamma_{n+1}(W_0)\simeq
\Gamma_{n+1}(W_1)$, and we can compare the $(\Gamma_{n+1},\Gamma_n)$ chain invariants.
If these agree then the $\Gamma_{n+1}$ chains are the same, and we have the data needed
to continue to the next stage.   It is these invariants that are related to the
Cochran-Orr-Teichner filtration of link concordance \cite{COT}.

The paper is organized as follows: Section 2 gives $\Gamma$ versions of the
chain and spine realization theorems of Part I. Section 3 describes in detail the move
used in deformations. This is a refinement of the move used in Part I. There is
also a generalization using capped gropes which is the geometric move appropriate
for solvable extensions.
 Section 4  defines
 $(\Lambda,\Gamma)$ chain invariants and suggests a possible reformulation using algebraic $K$-theory of
noncommutative localizations. In Section 5 we state relative versions of the
realization theorems, and give proofs. 
In Section 6 the general theory is specialized to the case of abelian
extensions, then iterated to give the obstruction theory for solvable extensions.

\head 2. Realization\endhead
In this section $N$ is a closed smooth 4-manifold. Versions in which $N$ may have boundary are given
in Section~4. Fix a homomorphism $\Gamma\to
\pi_1N$. A $\Gamma$ decomposition of $N$ is a dual decomposition $N=M\cup W$ together with a factorization
$\pi_1W\to \Gamma\to \pi_1N$, so that $\pi_1W\to \Gamma$ is onto. Note that for a $\Gamma$ decomposition to
exist $\Gamma$ must be finitely generated and $\Gamma\to \pi_1N$ must be onto, so we assume these
conditions hold as well.
\proclaim{2.1 Theorem ($\Gamma$-Chain Realization)} Suppose $D_*$ is a finitely generated based free chain
complex over $Z[\Gamma]$  with a chain map $f\:D_*\otimes Z[\pi_1N]\to C_*(N;Z[\pi_1N])$. Then
there is a $\Gamma$ decomposition $(N = M\cup W,\pi_1W\to \Gamma)$ realizing $(D,f)$ if and only
if
$D$ is homologically 2-dimensional and $H_0(D;Z[\Gamma])\to H_0(N;Z[\pi_1N])$
$ (=Z)$ is an isomorphism.\endproclaim
``Realizing $(D,f)$''  means there is a simple equivalence $g\:D_*\to
C_*(W;Z[\Gamma])$ so that the diagram
$$\CD D_*\otimes Z[\pi_1N]@>f>> C_*(N;Z[\pi_1N])\\
@VV {g} V @VV=V\\
C_*(W;Z[\pi_1N])@>{\text{inclusion}}>> C_*(N;Z[\pi_1N])\endCD$$
commutes up to chain homotopy. Other terms, e\.g\. ``homologically 2-dimensional'' are defined in Part~I.

\subhead 2.2 $\Gamma$-Spine Realization \endsubhead
Theorem 2.1 describes the chain complexes that can be realized. In the next statement the
chains are held fixed and
 the spine is varied using $\Gamma$ deformations. Deformations are described in detail in Section
3; for the next statement we need the following:
\roster\item $\Gamma$ deformations of $M$ are defined for codimension-0 submanifolds $M\subset N$ with a
factorization
$\pi_1(N-M)\to
\Gamma\to \pi_1N$. In particular if 
$N=M\cup W$,
$\pi_1W\to
\Gamma\to \pi_1N$ are a $\Gamma$ decomposition in the sense of the introduction then $M$ can be $\Gamma$
deformed. 
\item The output is another $\Gamma$ decomposition $M'\cup W'$, $\pi_1W'\to \Gamma\to
\pi_1N$.
\item There is a canonical simple chain equivalence $C_*(W;Z[\Gamma])@>\sim>>C_*(W';Z[\Gamma])$, whose
$Z[\pi_1N]$ reduction chain homotopy commutes with the inclusions to $C_*(N;Z[\pi_1N])$.
\endroster
We show that the $W$ part of the decomposition can change fairly arbitrarily within these constraints. What
makes this interesting is that the  changes in
$M$ are  more subtle: the spine is unchanged up to homotopy 2-deformation, for instance.

\proclaim{$\Gamma$-Spine Realization Theorem} Suppose $N=M\cup W$ is a dual decomposition, $K\to N$
is a 2-complex, and factorizations $\pi_1W\to \Gamma\to \pi_1N$ and $\pi_1K\to \Gamma\to \pi_1N$ are given.
 Then there is a $\Gamma$ deformation of $M$ to a decomposition $M'\cup W'$ so that the spine of\/ $W'$
realizes $K$ up to homotopy deformation, if and only if there is a simple chain
equivalence
$C_*(K;Z[\Gamma])\to C_*(W;Z[\Gamma])$ whose $Z[\pi_1N]$ reduction chain-homotopy commutes with the
inclusions.
\endproclaim
In fact the simple chain equivalence is realized in the sense that the diagram
$$\CD C_*(K;Z[\Gamma])@>\text{given}>> C_*(W;Z[\Gamma])\\
@VV{=}V @VV{\Gamma \text{deformation}}V\\
C_*(K;Z[\Gamma])@>{\text{homotopy deformation}}>> C_*(W';Z[\Gamma])\endCD$$
chain homotopy commutes.

Recall from Part I that a homotopy deformation of 2-complexes is a sequence of moves, either expansions or
collapses of 1- or 2-cells, or changes of attaching map of 2-cells by homotopy. A homotopy deformation
determines a simple homotopy equivalence, which in turn induces a simple chain
equivalence (with any coefficients). The converses are not true: a simple chain equivalence between
cellular chains of spaces is rarely induced by a map, and we do not expect simple homotopy equivalences of
2-complexes  to always come from homotopy 2-deformations. This last statement is the ``Andrews-Curtis
conjecture'', and  though expected to be false it is still open.

\subhead 2.3 $\Gamma_n$-Spine Realization \endsubhead
Grope deformations are  refinements that geometrically encode some  fundamental group
information.  They therefore apply only to certain
factorings of $\pi_1W$. Suppose
$\pi_1W\to
\Gamma$ as usual, and as in the introduction we define 
$$\Gamma_n(W) = (\pi_1W)/\text{ker}(\pi_1W\to\Gamma)^{(n)}.$$
Here $K^{(n)}$ is the $n$-fold iterated commutator subgroup defined by $K^{(0)}=K$ and $K^{(n)} =
[K^{(n-1)},K^{(n-1)}]$. Grope deformations are defined in \S3.5, and in \S3.6 we show that a $\Gamma$ grope
deformation of height
$n$ is a particular case of standard $\Gamma_{n-1}(W)$ deformation. The converse is not true, but  there
are enough of them to change skeleta in the same way.
\proclaim{Theorem} Suppose $N=M\cup W$ is a dual decomposition with factorization $\pi_1W\to \Gamma\to
\pi_1N$. Then the conclusions of Theorem 2.2 for the factorization $\pi_1W\to \Gamma_{n}W\to
\pi_1N$ hold with ``$\Gamma$ grope
deformation of height
$n+1$'' replacing ``$\Gamma_{n}W$-deformation''.
\endproclaim
The standard deformations are grope deformations of height 1, so this statement includes 2.2 as the
special case $n=0$.  As in 2.2 the deformations realize the chain equivalence.

\head 3.  $\Gamma$-deformations\endhead
In this section we suppose $N$ is a compact smooth 4-manifold with boundary divided into
pieces $\partial N = \partial_0M\cup \partial_0W$. At this point $M$ and $W$ are not defined; the notations
for the pieces of $\partial N$ are for later convenience. Fix a homomorphism
$\Gamma\to
\pi_1N$, and sometimes a further homomorphism $\Lambda\to \Gamma$. 
 We consider
extensions of the boundary decomposition to dual decompositions $N=M\cup W$, together with a factorization
$\pi_1W\to \Gamma\to \pi_1N$. Explicitly, ``extension'' of the boundary decomposition means $M\cap\partial
N=\partial_0M$, $W\cap\partial
N=\partial_0W$, and the significance of ``dual decomposition'' is that we are given handlebody structures
with handles of index $\leq 2$ on the pairs $(M,\partial_0M)$ and $(W,\partial_0W)$. The basic $\Gamma$ and $(\Lambda,\Gamma)$ deformation moves are described from several viewpoints in \S\S 3.1--3.3.
Standard structure (chain maps etc.) is derived in 3.4, and the grope version is described in~3.5--3.7.

A $\Gamma$ {\it deformation\/}  is a sequence of $\Gamma$ decompositions of $N$, each obtained from
the previous one by either ordinary handle moves in $M$ and $W$, or by the move described next:

\subhead 3.1 The data \endsubhead
As above we have $N=M\cup W$ and a factorization $\pi_1W\to
\Gamma\to \pi_1N$. The data for a  $\Gamma$-deformation of $M$ is:
\roster\item a compact orientable surface $\Sigma$ with boundary $S^1$ and a standard
hyperbolic basis of embedded curves;
\item an embedding $\Sigma\cup D^{\pm}_*\to \partial_1M^{(1)}$, where $D^{\pm}_*$ are
2-disks attached on the basis curves and
$\partial_1M^{(1)}$ is the level in the handlebody structure on $M$ between the
1-handles and the 2-handles; 
\item $\Sigma$ is disjoint from the 2-handles of $M$, so lies in $\partial_1W$.
This map $\Sigma \to W$ lifts to the $\Gamma$ cover of $W$; 
\item the disks are disjoint from the 2-handles of $W$; and
\item $\Sigma$ intersects a single 2-handle of $W$, and intersects the dual core of that
handle in a single point.
\endroster
If further a {\it pair\/} of groups $\Lambda\to
\Gamma\to\pi_1N$ is given and $\pi_1W$ lifts to $\Lambda$ then the data for a $(\Lambda,
\Gamma)$-deformation is the above and 
\roster\item[6] The boundary curve of $\Sigma$ lifts to the $\Lambda$ cover of $W$.
\endroster
This data is illustrated in Figure 1. In standard 4d terminology (cf.\ \cite{FQ})
$\Sigma\cup D^{\pm}_*$ is a disk-like capped surface. The caps intersect 2-handles of
$M$, while the body has a single point of intersection with (cocores of) handles of
$W$. Note that ``2-handles of $W$'' are by definition attached to lower handles of $W$.
The duals of these are handles attached to the complement (more specifically to
$\partial_1M^{(1)}$). The intersection with $\Sigma$ is with the attaching circle of
one of these dual handles.

\figure1

\subhead 3.2 The move as seen from $M$ \endsubhead
Suppose we have the data specified above, and denote by $A$ the handle of $W$
intersecting
$\Sigma$. We denote the dual handle by $A^*$, so it is the attaching circle of
$A^*$ that intersects $\Sigma$ in a point. The move changes $M\cup A^*$ by handle
moves and the result of the move ($M'$) is the complement of $A^*$ in
the resulting handlebody structure. The handle moves are obtained by pushing a
neighborhood of
$\Sigma\cup D^{\pm}_*$ over $A^*$. In 3.2.1 we describe this in a symmetric form, then in 3.2.2 a less
symmetric form that may be easier to see.

\subsubhead 3.2.1 The symmetric view\endsubsubhead A neighborhood of $(\Sigma\cup
D^{\pm}_*, \partial \Sigma)$ in $\partial_1M^{(1)}$ is isomorphic to a ball
$(D^2\times I, S^1\times I)$. The attaching circle of $A^*$ intersects this in an arc
$\{p\}\times I$. Push this ball across $A^*$ through the center to the other side,
dragging the attaching maps of
$M$ 2-handles along. After this the attaching maps miss the center of the core of
$A^*$, so we can deform them radially off $A^*$. This sequence of moves is shown in
Figure 2.

\figure2

This gives an isotopy of attaching maps of 2-handles of $M$, in the boundary of
$M^{(1)}\cup A^*$. Use this to do handle moves in the handlebody $M\cup A^*$. This
changes the handlebody structure but does not change the manifold. In particular it is
still embedded in $N$ with complement $W-A$. Define $M'$ to be the complement of $A^*$
in this new handlebody, and define $W'$ to be (the closure of) the complement of $M'$
in $N$. Note that $W'$ is given as $W-A$ union with the repositioned $A$.

\subsubhead 3.2.2 The asymmetric view\endsubsubhead We now give a less symmetric but
more elementary description of the deformation. Consider the top $D^2$ of the 
$D^2\times I$ neighborhood of the capped surface. This disk intersects the attaching
map of $A^*$ in a single point, and many attaching maps of $M$ 2-handles. Do handle
moves of the $M$ 2-handles across $A^*$ to make them disjoint from the disk. This
involves pushing along arcs in the disk, then over $A^*$. This operation is shown in
Figure~3. 

\figure3

Again we define $M'$ to be the resulting handlebody, with $A^*$ omitted.  Comparing
Figures 2 and 3 shows the two operations are isotopic, so this gives the same result
as before. Also 
we could do the same thing with the lower $D^2$ in the ball neighborhood. This also is
isotopic to the symmetric version, so we additionally see that using the upper or
lower disks yield isotopic results.

\subhead 3.3 The move as seen from $W$ \endsubhead
The surface $\Sigma$ intersects the core of the handle $A^*$ in a point, so intersects
$A^*$ in a disk. Deleting the interior of this disk leaves a surface $\Sigma_0$ with
two boundary circles. $\Sigma_0$ is disjoint from the interiors of all 2-handles so
it also lies in $\partial_1W^{(1)}$, the level between the 1- and 2-handles of $W$. In
this level the new boundary circle in $\Sigma_0$ is (a canonical parallel of) the
attaching circle of $A$. Duals of the 2-handles in $M$ are attached to this level, and
the disks $D^{\pm}_*$ go over these duals. 

Thicken $\Sigma_0\cup D^{\pm}_*$, in the upper boundary of $W\cup(\text{2-handles of
}M)$. This is the thickening as used in 3.2, minus a plug where $A$ passes through, so
is $(S^1\times I)\times I$. The top of this is an annulus with one end $\partial A$
and going over duals of handles of $M$. The asymetric description of the move given
in 3.2.2 is to push the handles of $M$ off this annulus across $A^*$. The dual of this
move is to move the attaching map of $A$ by isotopy along this annulus, going over
duals of $M$-handles in the process. In particular the $W'$ obtained at the end of the
move is constructed by removing the handle $A$ from $W$ and reattaching it on the other
end of the surface $\Sigma_0$. 

This description shows the link giving attaching maps in $W$ change by a variant of the ``grope cobordism''
of Conant and Teichner \cite{CT}. Here the grope has height 1 (is just a surface), and the standard basis
curves (attaching curves for caps) are required to be parallel to attaching maps of (duals of) handles of
$M$. Grope cobordisms of greater height appear in the grope moves in~3.5.

\subhead 3.4 Standard structure \endsubhead 
The ``standard structure'' of a deformation (see the beginning of 2.2) is 
\roster\item a canonical $\Gamma$  structure for the
output decomposition, i.e\. a factorization $\pi_1W' \to \Gamma \to \pi_1N$; and 
\item a simple chain equivalence of chain complexes
$$C_*(W,\partial_0W;Z[\Gamma])\to C_*(W',\partial_0W;Z[\Gamma]).$$
\endroster
A $(\Lambda,\Gamma)$ deformation has slightly more, namely
\roster\item[3] a canonical refinement of the $\Gamma$-decomposition structure to a $\Lambda$
structure, i.e\. the factorization above lifts canonically to $\pi_1W'\to \Lambda$.
\endroster
 We now derive this structure.

Recall that a $\Gamma$ deformation  is a sequence of decompositions, each obtained from the
previous one by either ordinary handle moves in $M$ and $W$, or by a move described above. We begin with the
ordinary handle  moves. These do not change $M$ and $W$, so $\pi_1W$ does not change and the original
$\pi_1$ factorization is used. The identity map $W\to W$ induces a chain map of complexes constructed using
the two handle decompositions. But it is a standard fact that this is a simple chain equivalence. Indeed the
definition of ``simple chain equivalence'' was developed exactly to encode the algebraic changes resulting
from handle moves, so this chain map is simple essentially by definition. 

It remains to describe structure for the new moves. On the $W$ side the move changes the attaching map of a
handle $A$. The new attaching map of
$A$ lifts to the $\Gamma$ cover, or $\Lambda$ cover in the $(\Lambda,\Gamma)$ case, so the homomorphism
$\pi_1(W-A)\to
\Gamma$ extends canonically to
$\pi_1W'
\to
\Gamma$ (respectively $\Lambda$). This provides the factoring.
Next, there is an evident bijection of handles in $W$ and $W'$. This bijection extends linearly to give
based isomorphisms of the cellular chain groups
$C_*(W,\partial_0W)\simeq C_*(W',\partial_0W)$. On the group  level this works with $\Lambda$
coefficients in the $(\Lambda,\Gamma)$ case. With $\Gamma$ coefficients this  is an isomorphism of chain
complexes, i.e\. commutes with the boundary homomorphisms. This is clear for all generators except $A$,
where it follows from the lifting of the surface $\Sigma_0$ as follows: Think of the chain groups as
relative homology of skeleta,
$$C_i(W,\partial_0W;Z[\Gamma])=H_i(W^{(i)},W^{(i-1)};Z[\Gamma]).$$
Then 
$\partial A$ is the image of the
attaching map in  $H_1(W^{(1)},W^{(0)};Z[\Gamma])$. But the lift of
$\Sigma_0$ to the $\Gamma$ cover provides a homology between the images of the two
attaching maps, so they are equal in this group.

\subhead 3.5 Grope moves\endsubhead
We first recall the definition of a grope (see \cite{FQ}). A grope of {\it height 1\/} is a capped surface,
as described above. It has a standard embedding in $D^3$ with the boundary of the surface a standard circle
in $S^2$. We now proceed recursively: suppose gropes of height $n-1$ are defined, and define gropes of
height $n$ to be the result of replacing all the caps in a grope of height $n-1$ by capped surfaces. The
simplest grope of height 2 is shown in Figure~4.

The {\it caps\/} of the grope are the caps of the layer of capped surfaces. The rest of the grope is called
the {\it body\/}, and consists of $n$ layers of surfaces. There is a standard model (embedding in $D^3$)
obtained by starting with the standard model of the grope of height $n-1$ and replacing $D^3$ neighborhoods
of its caps by standard models of capped surfaces. Further these are all spines of $D^3$, or equivalently,
the regular neighborhood of the grope is isomorphic to $D^3$. Embeddings and immersions of gropes are
defined in terms of these standard neighborhoods: an embedding of a grope in a 4-manifold, for instance, is
an embedding of $D^3\times I$, thought of as a neighborhood of the spine.
\subsubhead 3.5.1 Data for a $\Gamma$ grope move\endsubsubhead
The data is the straightforward analog of 3.1:
\roster\item a capped grope of height $n$ embedded in $\partial_1M^{(1)}$;
\item the body of the grope is disjoint from the attaching maps of the 2-handles of $M$, so lies in
$\partial_1W$. This inclusion lifts to an inclusion of the body in the $\Gamma$ cover;
\item the cap disks are disjoint from the 2-handles of $W$; and 
\item the body intersects (dual attaching maps of) 2-handles of $W$ in a single point in the lowest level
surface.
\endroster
This is illustrated in Figure 4. 

\figure4

\subsubhead 3.5.2 The move\endsubsubhead
The move itself is the straightforward elaboration of the height-1 case described in 3.2 and 3.3. Denote by
$A^*$ the dual of the handle of $W$ that intersects the grope body, and arrange that its attaching map
intersects the neighborhood
$D^3\times I$ in an arc $\{x\}\times I$. The move pushes the $D^3\times I$ and all the stuff intersecting
it across $A^*$, as in 3.2. Alternatively, it relocates the handle $A\subset W$ to one attached on the
boundary circle of the grope body, as in~3.3.

\subhead 3.6  Contraction, and standard structure for grope moves\endsubhead
The (one-fold) {\it contraction\/} of a grope undoes the last stage in the construction, which is to say
replaces the uppermost layer of capped surfaces by disks. This can be described explicitly in terms of the
grope, by using some of the caps to do surgery on the uppermost surfaces, reducing them to disks. The
$k$-fold contraction does this $k$ times. Thus the $k$-fold contraction of a grope of height $n$ gives a
grope of height $n-k$, whose body is the first $n-k$ layers of surfaces in the original.

Contracting the grope in the data for a grope move gives data for a lower-height grope move. Geometrically
the two moves have the same effect since the two gropes have the same regular neighborhood. They agree
algebraically as well after taking into account a
canonical improvement in the fundamental group structure. If  $\pi_1W$ factors through $\Gamma$ then
$\Gamma_k(W)$ is defined (introduction and 2.3) to be
$\pi_1W$ modulo the $k$-fold iterated commutator of the kernel of $\pi_1W\to \Gamma$. 

\proclaim{Proposition} Suppose the body of a grope of height $n$ lifts to the $\Gamma$ cover of $W$. Then
the union of the first $n-k$ layers lifts to the $\Gamma_k(W)$ cover.\endproclaim
 \demo{Proof}Suppose we have an oriented
surface with one boundary component. Then in $\pi_1$ the boundary curve is a product of commutators. More
generally, the boundary curve of a grope of height $k$ is an $k$-fold iterated commutator in the
fundamental group of the body. Now suppose $G$ is a
capped grope of height
$n$ with a map of the body to $W$ that lifts to the $\Gamma$ cover. This means the attaching curves for
the caps lie in the kernel of $\pi_1W\to \Gamma$. The attaching curves on the body of the
$k$-fold contraction are $k$-fold commutators of the original cap curves, so  this body  lifts to the
$\Gamma_k(W)$ cover of $W$. \enddemo

Setting $k=n-1$, and recalling that the that the standard moves are grope deformations of height 1 we get:
\proclaim{Corollary} A\/ $\Gamma$ grope deformation of height $n$ is a 
$\Gamma_{n-1}(W)$-deformation.\endproclaim

There is a weak converse to this: suppose a grope body of height $n-k$ lifts to the $\Gamma_{k}(W)$ cover.
Then it extends to a {\it map\/} of a grope body of height $n$ that lifts to the $\Gamma$ cover. In general
this map will not be suitable for a deformation move because it is not an embedding. It may be possible to modify it to get an embedding, cf\. I\S7 and 5.4 below. 

\subhead 3.7 Relation to finite type filtrations\endsubhead
``Finite type'' invariants of links are defined to be ones that do not detect changes of sufficient
complexity. Specifically a ``degree $k$''  modification
is defined as follows: suppose a 3-ball intersects the link in $k$ segments. Change the intersection with the ball
to a new configuration with the property that if any one segment is omitted then the result is isotopic (in
the disk, rel boundary) to the corresponding omission from the original. An invariant is said to be ``finite type of degree $<k$'' if it is unchanged by degree-$k$ modifications. More precisely an invariant defined for a class of links is finite type of degree $k$ if the class is closed under degree $k+1$ modifications and the invariant is unchanged by these modifications. 

``Finite type'' is related to the grope picture as follows. The simplest capped surface is a once-punctured
torus with a single pair of caps. Similarly the simplest grope of height
$n$ is obtained by always replacing caps with copies of this simplest capped surface. Figure 4 shows the
simplest grope of height 2. Arbitrary grope moves can be subdivided into a sequence of moves all using the
simplest grope, so we focus on these. Such a grope has $2^n$ caps. We can further suppose that the
embeddings used as data for a move have each cap intersecting a single $M$ 2-handle in a point. The
neighborhood of the grope is then a disk that intersects the attaching link of the 2-handles in $2^n$
segments. The grope move changes the attaching link by pushing it across $A^*$. 
But in fact we could pull all but one of these out of the disk, to arrange the intersection to be
$2^n$ sub-segments of one of the original segments. After this the grope move changes only one of the
original segments. Thus omitting any one of these segments trivializes the move. The conclusion is that
grope moves of height $n$ on a link are modifications of degree $2^n$.  Conant and
Teichner \cite{CT} have developed a sharper relation between gropes and finite type that  characterizes modifications of arbitrary degree. 

 This discussion shows that the solvable-tower version of chain invariants described in \S6 are in an appropriate sense finite type invariants.

\head 4. Chain invariants\endhead
In this section we define an invariant using equivalence classes of algebraic data coming from a
$(\Lambda,\Gamma)$ deformation. This should be regarded as preliminary: even if it gives the right thing
the formulation needs refinement. See 4.6 for a possible relation to $K$-theory.

\subhead 4.1 Small chain objects and $\sc(Z[\Lambda])\downarrow C$\endsubhead
 We will say $D$ is a ``small chain object over $Z[\Lambda]$'' if it satisfies the
following:
\roster\item $D$ is a finitely generated free based
$Z[\Lambda]$ chain complex;
\item $D$ is  homologically 2-dimensional; and
\item $H_0(D;Z[\Lambda])$ is trivial as a $Z[\Lambda]$ module.
\endroster
Small chain objects encode the algebraic properties of chains of a 2-complex, or relative chains of a
relative 2-complex. In particular if $M\cup W$ is a $\Lambda$ decomposition then the cellular chain complex
$C_*(W,\partial_0W;Z[\Lambda])$ is a small chain object.

Denote by $\sc(Z[\Lambda])$ the category of small chain objects and simple chain
equivalences. A homomorphism $\Lambda\to \Gamma$ induces a functor of categories by
$D\mapsto D\otimes Z[\Gamma]$. 

Now suppose $C$ is a small chain object over $Z[\Gamma]$. Then we denote the fiber over $C$ of the functor
$\sc(Z[\Lambda])\to \sc(Z[\Gamma])$  by $\sc(Z[\Lambda])\downarrow
C$. Explicitly this is the category with objects $(D,f)$ where
$D$ is a small chain object over
$Z[\Lambda]$ and $f$ is a simple chain equivalence $D\otimes Z[\Gamma] \to C$. Morphisms
are simple chain equivalences over $Z[\Lambda]$ whose $Z[\Gamma]$ reductions chain homotopy commute with
the given equivalences. 

The morphisms define an equivalence relation on objects of $\sc(Z[\Lambda])\downarrow C$.
We abuse the notation  by using $\sc(Z[\Lambda])\downarrow C$ to  refer also to the
set of equivalence classes. The basic plan is to use the equivalence class of the chain complex of $W$ as our invariant. The definition above is sufficient when $\partial_0W$ is empty, but has to be elaborated a bit in general. This is done in the next two sections. 

\subhead 4.2 Chains with stratified coefficients\endsubhead
When $\partial_0W$ is nonempty we need
{\it stratified\/} coefficients to record $\pi_1$ information on the boundary. In fact we use the language
``stratified coefficients'' but only describe the very special case needed
rather than develop a general theory and specialize it.

The ``suspension'' of a complex is defined by shifting the groups down one degree and multiplying the
boundary homomorphisms by $(-1)^n$. In detail the suspension of $C_*$ is denoted by $C_{*-1}$. The
$n^{\text{th}}$ chain group of  $C_{*-1}$ is $C_{n-1}$, and the boundary from degree $n$ to degree $n-1$ is
$(-1)^{n-1}\partial\:C_{n-1}\to C_{n-2}$. 

If $(X,Y)$ is a CW pair then there is a boundary chain map of cellular chain complexes
$d\:C_*(X,Y)\to C_{*-1}(Y)$ that on homology induces the boundary homomorphism in the long exact
sequence. Here the target is the suspension as defined above, and the signs on the boundary homomorphisms
are needed to make this a chain map.

Now suppose
$\lambda\:\Lambda\to \Gamma$ is a homomorphism of groups, $(X, Y)$ is a CW pair, and
there is a commutative diagram
$$\CD\pi_1Y@>>>\Lambda\\
@VVV@VV{\lambda}V\\
\pi_1X@>>>\Gamma.\endCD$$
Then we define the stratified chain complex $C_*(X,Y;Z[\Lambda])$ to be the pullback of the diagram
$$\CD @. C_{*-1}(Y;Z[\Lambda])\\
@.@VV{\lambda}V\\
C_*(X,Y;Z[\Gamma])@>d >> C_{*-1}(Y;Z[\Gamma])\endCD$$
Explicitly the $n^{\text{th}}$ chain group in the pullback is given by $(x,y)\in C_n(X,Y;Z[\Gamma])\oplus
C_{n-1}(Y;Z[\Lambda])$ such that $d(x)= \lambda(y)$, where $\lambda$ denotes change-of-coefficient map on
chains induced by the homomorphism $\lambda$. 

The facts we need about this complex are:
\roster\item $C_{*}(X,Y;Z[\lambda])$ is a chain complex over $Z[\Lambda]$, though the chain groups are not
free;
\item $C_{*}(X,Y;Z[\lambda])\otimes Z[\Gamma]$ is $C_{*}(X,Y;Z[\Gamma])$, so is free and based over
$Z[\Gamma]$; and
\item if the homomorphism $\pi_1X\to \Gamma$ lifts to $\Lambda$ compatibly with the homomorphism on
$\pi_1Y$, then there is a natural $Z[\Lambda]$ chain map  $C_{*}(X,Y;Z[\Lambda])\to C_{*}(X,Y;Z[\lambda])$
induced by  the universal property of the pullback and the commutative diagram 
\endroster
$$\CD C_*(X,Y;Z[\Lambda]) @>d>> C_{*-1}(Y;Z[\Lambda])\\
@VV{\lambda}V @VV{\lambda}V\\
C_*(X,Y;Z[\Gamma])@>d >> C_{*-1}(Y;Z[\Gamma])\endCD$$

\subhead 4.3 Small chain objects over  stratified chains \endsubhead
The definition of the category $\sc(Z[\Lambda])\downarrow
C$ given in 4.1 requires $C$ to be a small chain object over $Z[\Gamma]$. Here we extend this slightly to
include stratified chains. Suppose $C$ is
\roster\item a chain complex (not necessarily free) over $Z[\Lambda]$, and
\item $C\otimes Z[\Gamma]$ is a small chain object over $Z[\Gamma]$. 
\endroster
Then we define $\sc(Z[\Lambda])\downarrow C$ to be the category with objects $(E,f)$, where $E$ is a small
chain object over
$Z[\Lambda]$ and $f\:E\to C$ is a $Z[\Lambda]$ chain map  whose $Z[\Gamma]$ reduction is a simple
equivalence. Morphisms in the category are simple equivalences over $Z[\Lambda]$ that chain homotopy
commute with the maps to $C$. 

As before we use the same notation for the set of equivalence classes of objects.

\proclaim{4.4 Proposition (Chain invariant realization)} Suppose $\Lambda@>\lambda>>\Gamma\to \pi_1N$ are given and
$N=M\cup W$ is a $\Lambda$-decomposition. Then the stratified
chains
$C=C_*(W,\partial_0W;Z[\lambda])$ are defined. If $f\:M\cup W\mapsto M'\cup W'$ is a $\Gamma$ deformation so
that the lift on $\pi_1\partial_0W$ extends to $\pi_1W'\to \Lambda$ then
$(C_*(W',\partial_0W;Z[\Lambda]),f_*)$ represents a class in $\sc(Z[\Lambda])\downarrow C$. Conversely any
such class is realized by some  deformation and lift.
\endproclaim
 
\proclaim{Refinement} If\/ $\Lambda\to\Gamma$ factors through\/ $\Gamma_{n}(W)$ then invariants can be
realized by $\Gamma$ grope deformation of height $n+1$.\endproclaim
\demo{Proof}  Putting together the $H_0$ hypothesis on small chain objects, the chain
equivalence, and the fact that
$M\cup W$ is a dual decomposition, we see the natural maps induce isomorphisms
$$H_0(D;Z[\Lambda])\simeq H_0(D;Z[\Gamma])\simeq H_0(W,\partial_0W;Z[\Gamma])\simeq
H_0(N,\partial_0W;Z[\pi_1N]).$$
This means, according to the Chain Realization Theorem for $\Lambda$, there is a
$\Lambda$-decomposition $\widehat W\cup \widehat M$ with a  simple chain equivalence $D_*@>\sim>> C_*(\widehat
W,\partial_0W;Z[\Lambda])$.

Let $(K,\partial_0W)\to (N,\partial_0W)$ be the spine of $\widehat W$. Then $K$ and the chain
equivalence over $\Gamma$ satisfy the hypotheses of the Spine Realization theorem. Thus
there is a $\Gamma$ deformation of $M\cup W$ to a decomposition $M'\cup W'$ and a homotopy
deformation of the  the spine of $W'$ to $K$. Further these realize the input data in the
sense that the canonical chain equivalence from the $\Gamma$ deformation is chain
homotopic to the input equivalence. We conclude that the decomposition $M'\cup W'$ and
its associated data does realize the given chain invariant.

The refinement is obtained by subsituting 5.3 for 5.2 in this argument.
\enddemo
\subhead 4.5 $(\Lambda,\Gamma)$ deformation\endsubhead
The $\pi_1W'$ lift to $\Lambda$ needed to define a chain invariant comes automatically in a
$(\Lambda,\Gamma)$ deformation. 
A natural question is therefore: are all chain invariants realized by $(\Lambda,\Gamma)$ deformation?
Proposition 4.4 does not show this, and we suspect  that there is
significant structure not identified in 4.3 that  limits the classes realized. For instance the
$(\Lambda,\Gamma)$ deformation moves come with a canonical basis-preserving isomorphism of the $Z[\Lambda]$
chain groups. This is usually not a chain map, but it may be close enough to one to show that some
variant of Whitehead torsion vanishes. A possible such variant is described below. Or a $(\Lambda,\Gamma)$
deformation may have a ``trace'' embedded in $N\times I$ that gives a concordance between the two
decompositions. In that case a bilinear invariant would be defined, and the linear invariants would be
limited to the image of the bilinear ones.

\subhead 4.6 $K$-theory\endsubhead
Suppose $\Lambda\to \Gamma$ is a homomorphism and $f\:D\to E$ is a chain map of
$Z[\Lambda]$ complexes that becomes a simple chain equivalence when tensored with
$Z[\Gamma]$. Then it is also a chain equivalence over the Cohn localization
$Z[\Lambda]_{\Sigma}$, obtained by inverting the set $\Sigma$ of square matrices over
$Z[\Lambda]$ whose reductions to
$Z[\Gamma]$ are isomorphisms. See Ranicki \cite{R} for an overview of topological applications of Cohn
localizations. The chain equivalence has a torsion
$\tau(f)\in K_1(Z[\Lambda]_{\Sigma})/\pm \Lambda$.  Here, as usual in the Whitehead group, we 
divide out the units $\pm \Lambda$ to make the torsion independent of choice of basis cells in the
$\Lambda$ cover. Since $f$ is simple over $Z[\Gamma]$ this torsion lies in the
kernel of the homomorphism to $K_1(Z[\Gamma])/\pm \Gamma$.

\subsubhead 4.6.1 Question\endsubsubhead Does a $(\Lambda,\Gamma)$ deformation define a torsion in 
$$\bigl(\text{ker}(K_1(Z[\Lambda]_{\Sigma})\to
K_1(Z[\Gamma])\bigr)/\pm \Lambda?$$

For the next question we suppose $h\:C\to D$ is a
$Z[\Lambda]$ chain map, $C$ is a small chain object over $Z[\Lambda]$, $D\otimes Z[\Gamma]$ is a small
chain object, and $h\otimes Z[\Gamma]$ is a simple equivalence. We are thinking of the $Z[\Lambda]$ chains
mapping to the stratified chains of a pair $(W,\partial_0W)$. Then
there is a natural map $\Theta_h$,
$$ \bigl(\text{ker}(K_1(Z[\Lambda]_{\Sigma})\to
K_1(Z[\Gamma])\bigr)/\pm \Lambda @>\Theta_h>> \sc(Z[\Lambda])\downarrow D.$$
defined by composing $\partial\:C_2\to C_1$ with a homomorphism $C_2\to C_2$ representing  the
$K$-theory class.  We expect that in general the chain invariant group  depends on $C$ and $D$ and not
just $\Lambda\to \Gamma$. However there may be useful cases in which they agree:

\subsubhead 4.6.2 Question\endsubsubhead
 Is $\Theta_h$  an isomorphism if either
\roster \item If $C_i=0$ for $i\neq 2$, or
\item $\Lambda=\Gamma_n(W)$ and $C\to D$ are the $\Lambda$ and stratified chains for some $\Gamma$
decomposition
$M\cup W$.
\endroster
The difference between the two groups concerns whether simple chain equivalences over $Z[\Gamma]$ can be
lifted in a reasonably canonical way to $Z[\Lambda]$ chain maps. The idea in the first case is that the
model is so small that not too much can go wrong. The reason this case is of interest is that it may include
decompositions manufactured to study more general links. In the second case the idea is that behavior in the chain complexes in degrees less than 2 is determined by
the group, so again problems are concentrated in degree 2. This case is interesting because it applies to
the solvable tower analysis.

\head 5. Relative versions, and proofs of realization\endhead
In this section we state relative versions of Theorems 2.1--2.3 and give  proofs. The proofs are mostly
minor modifications of the proofs of Part~I.

As in Section 3 we fix a compact smooth 4-manifold $N$ and a decomposition of its boundary into $\partial N
=
\partial_0M\cup
\partial_0W$. Again, $M$ and $W$ are not yet defined; the notations
for the pieces of $\partial N$ are for later convenience. We also fix a homomorphism
$\gamma\:\Gamma\to\pi_1N$ and a factoring $\pi_1\partial_0W\to \Gamma\to\pi_1 N$,
 and consider $\Gamma$-decompositions $N=M\cup W$ that extend the structure on
$\partial N$. 

The first topic is the extension of the chain realization theorem 2.1. Note that the factoring of
$\pi_1\partial_0W$ through $\gamma\:\Gamma\to\pi_1N$ means the stratified chains
$C_*(N,\partial_0W;Z[\gamma])$ are defined. Recall that a ``small chain object'' is a free based
homologically 2-dimensional complex with $H_0$ trivial as a $Z[\Gamma]$ module.

\proclaim{5.1 Theorem (Relative $\Gamma$-Chain Realization)} A $\Gamma$ decomposition $N=M\cup W$ defines a
small chain object $D=C_*(W,\partial_0W;Z[\Gamma])$ over
$Z[\Gamma]$ and a $Z[\Gamma]$ chain map $f\:D\to C_*(N,\partial_0W;Z[\Gamma])$  whose $Z[\pi_1N]$
reduction is an isomorphism on $H_0$ and an epimorphism on $H_1$. Conversely given such $D$ and $f$ 
there is a
$\Gamma$ decomposition  that realizes it up to simple chain homotopy.
\endproclaim

The first step is to modify the 1-skeleton alignment lemma of Part I \S4. The modifications are relatively
minor, but the statement is complicated enough that to be safe we repeat the whole thing.
\proclaim{5.2 1-skeleton alignment Lemma} 
Suppose $\Gamma\to\pi$,
 $C_*$ is a free based $Z[\Gamma]$ complex whose $H_0$ is trivial as a $Z[\Gamma]$ module, $D_*$ is a free
based
$Z[\pi]$ complex, and
$f\:C_*\to D_*$ is a
$Z[\Gamma]$ chain map whose $Z[\pi]$ reduction is an isomorphism on $H_0$ and an epimorphism on $H_1$. Then
there is a chain homotopy commutative diagram
$$\CD C_*@>f>>D_*\\
@VVV@VVV\\C'_*@>f'>> D'_*\endCD$$
so the vertical maps are simple equivalences over $Z[\Gamma]$ and $Z[\pi]$ respectively and are
isomorphisms in degrees $\geq3$, and $f'$ is basis-preserving in degrees 0 and 1. Further:

If  $D_*$ is the cellular chains of a CW complex or pair (resp.\ 4-d handlebody) with connected $\Gamma$
cover then $D_*\to D'_*$ can be arranged to be induced by a homotopy 2-deformation (resp.\ handlebody
moves).

If both $C_*$ and $D_*$ are cellular chains of CW complexes or pairs (resp.\ 4-d handlebodies) with
connected covers and the isomorphism on $H_0$ is the identity, then both $C_*\to C'_*$ and $D_*\to D_*$ can
be arranged to be induced by 2-deformation (resp.\ handle moves).
\endproclaim
To make sense of the requirement that $f$ be a $Z[\Gamma]$ chain map we regard $D_*$ as a $Z[\Gamma]$
complex via the module structures induced by the ring map $Z[\Gamma]\to Z[\pi]$. Such a chain map  factors
through the quotient
$C_*\to C_*\otimes Z[\pi]$ and a
$Z[\pi]$ chain map
$C_*\otimes Z[\pi]\to D_*$. This latter chain map is the ``$Z[\pi]$ reduction'' of $f$.

\demo{Proof}  The proof is obtained by modifying the proof of  I\S4.
 In the geometric cases the new hypothesis are automatically satisfied,
and the new conclusion follows for free  because  cell arguments lift to any cover.  The cases that
require modification are
 ``semi-algebraic 0-skeleton alignment'' and ``algebraic 1-skeleton alignment.'' 

In the semi-algebraic 0-skeleton argument  a copy of $D_0$ is added to $C_0$. In order to get a free
$\Gamma$ complex we use a free $Z[\Gamma]$ module $\hat D_0$ with the same basis as $D_0$. Next we must
lift $f_0$ and $g_0$ to homomorphisms $\hat f_0$ and $\hat g_0$ with $Z[\Gamma]$ coefficients, so that
 equation (2) 
$$\hat f_0 \hat g_0 +\partial t = \id_C$$
holds with $Z[\Gamma]$ coefficients. Previously it was done with $Z[\pi]$ coefficients. To get the
improvement we need that the change-of-coefficient homomorphism $C\to C\otimes Z[\pi]$ induces an
isomorphism on $H_0$. But this follows from the  hypothesis that $H_0$ is trivial as a $Z[\Gamma]$ module.
After this adjustment the argument proceeds as before.

In the algebraic 1-skeleton alignment we take $\hat D_1$ a free $Z[\Gamma]$ module whose $Z[\pi]$ reduction
is $D_1$. Lift $g\:D_1\to C_1\otimes Z[\pi]$ to $\hat g\:\hat D_1\to C_1$. Then the right two columns in
the diagram below equation (1)  become simple equivalences of $Z[\Gamma]$ complexes. The
lower vertical arrow in the second column becomes $(\partial, \partial\hat g)$. This does not affect
the rest of the argument since the $Z[\pi]$ reduction of $\partial\hat g$  can be identified with
$\partial$, as was done in the original argument.

This ends the modification of the proof of the  1-skeleton alignment lemma.
\enddemo
\demo{Proof of 5.1}
The rest of the proof of  5.1 is obtained by routine modifications of Part I \S5. These are even more minor
than the modifications needed in 1-skeleton alignment, so are omitted.
\enddemo

Next is the bounded version of spine realization. Again we describe data coming from a deformation, and the
theorem asserts that any such data is realized. The main difference is the appearence of  stratified
chains encoding the fixed structure on the boundary. Suppose $d\:M\cup W\mapsto M'\cup W'$ is a $\Gamma$
deformation, and denote the (relative) spine of $W'$ by $(K,\partial_0W)$. Then
\roster\item  $(K,\partial_0W)\to N$ is a relative 2-complex and the map is the identity on $\partial_0W$;
\item there is a factoring of $\pi_1K$ through $\Gamma$ that extends the given factoring on
$\pi_1\partial_0W$; and
\item the deformation induces a simple  chain equivalence $d_*\:C_*(K,\partial_0W;Z[\Gamma])\to
C_*(W,\partial_0W;Z[\Gamma])$ that chain homotopy commutes with the maps to  the stratified chains
$C_*(N,\partial_0W;Z[\gamma])$.
\endroster

\proclaim{5.3 Theorem (Relative $\Gamma$-Spine Realization)} Suppose $N=M\cup W$ is a $\Gamma$
decomposition,
$(K,\partial_0W)\to (N,\partial_0W)$ is a relative 2-complex, and $$d_*\:C_*(K,\partial_0W;Z[\Gamma])\to
C_*(W,\partial_0W;Z[\Gamma])$$ is a chain map satisfying the conditions above. 
 Then there is a $\Gamma$ deformation of $M$ to a decomposition $M'\cup W'$ so that the spine of\/
$(W',\partial_0W)$ realizes $(K,\partial_0W)$ up to homotopy deformation.
\endproclaim
As in 2.2  the homotopy deformation can be arranged to realize the  chain data.

{\it Proof\/}
The proof closely follows the proof in Part I \S7. Note, however, that $M$ and $W$ are interchanged in the
notation used here and in I~\S7.
\subsubhead{5.3.1 Align 1-skeleta}\endsubsubhead Use the improved 1-skeleton alignment lemma above 
to change the handle structure on $(W,\partial_0W)$ and homotopy deform $(K,\partial_0W)$ so the 1-skeleta
can be identified. Specifically we get
\roster\item $K=W^{(1)}\cup(\text{2-cells)}$;
\item the homomorphisms $\pi_1K\to \Gamma$ and $\pi_1W\to \Gamma$ agree on the 1-skeleton;
\item the simple chain equivalence $d_*\:C_*(K,\partial_0W;Z[\Gamma])\to C_*(W,\partial_0W;Z[\Gamma])$ is
the identity on $C_0$ and $C_1$, and is a basis-preserving isomorphism on $C_2$; and
\item the chain maps to the stratified chains $C_*(N\partial_0W;Z[\gamma])$ agree on the 0- and 1-chains,
and commute with $d_*$.
\endroster

Note that $\Gamma$-deformations (of
$M$) preserve all this data. Deformations do not change the 1-skeleton of $W$, so do not affect
(1) and (2). The 2-handles do change, but the standard data gives a basis-preserving chain map of
$Z[\Gamma]$ chains from the old to the new, so composing this with the chain equivalence (3) replicates
this chain data. Finally the standard data chain map is the identity on the 1-skeleton, so condition (4) is
not disturbed either.

\subsubhead{5.3.2   2-cell data}\endsubsubhead
To work with 2-cells we need the absolute rather than relative chain complexes. It is  this step that
requires use of stratified chains. The relative chains of both $K$ and $W$ map to the stratified chains of
$(N,\partial_0W)$, and this in turn has a boundary map to the $Z[\Gamma]$ chains of $\partial_0W$. This
defines the right-hand square in the diagram
$$\CD C_*(K)@>>>C_*(K,\partial_0W)@>>>C_{*-1}(\partial_0W)\\
@VV{\hat d_*}V @VV{d_*}V @VV=V\\
C_*(W)@>>>C_*(W,\partial_0W)@>>>C_{*-1}(\partial_0W)\endCD$$
where all coefficients are $Z[\Gamma]$. The  right-hand square commutes, so induces
the left vertical chain map $\hat f_*$. This also is the identity on the 1-skeleton and a basis-preserving
isomorphism on the 2-cells.

The  chain map $\hat d_*$ gives  $Z[\Gamma]$ homologies between  attaching
maps of the 2-handles of $W$ and the 2-cells of $K$. To see these recall that the cellular chain complex is
defined to have chain groups $C_i(K)=H_i(K^{(i)},K^{(i-1)})$, where $K^{(i)}$ denotes the $i$-skeleton, and
boundaries are defined by  compositions
$$H_i(K^{(i)},K^{(i-1)})@>\partial>>H_{i-1}(K^{(i-1)})@>>>H_{i-1}(K^{(i-1)},K^{(i-2)}).$$
Since $d_2$ takes 2-cells to 2-handles, their boundaries must be equal in $H_1$ of the 1-skeleton.

We claim that to prove the theorem it is sufficient to deform $W$ so these 
$Z[\Gamma]$ homologies
are realized by homotopies. First, the homotopies give the desired 2-deformation between $K$ and the spine
of $W$. However we also want this 2-deformation to homotopy commute with the maps to $N$, and for this  we
need to extend the homotopies of attaching maps to homotopies between the maps of the 2-cells into $N$. The
obstruction to finding such a homotopy is an element of $\pi_2N$ formed from a 2-cell in each of $K$ and
$W$, and the homotopy between the attaching maps. Since we are working in the universal cover of
$N$,
$\pi_2N=H_2(N;Z[\pi_1N])$ and it is sufficient to show the homology class vanishes. 

The inclusion $(W,\partial_0W)\to (N,\partial_0W)$ induces a commutative diagram 
$$\CD C_*(W;Z[\pi_1N])@>>>C_*(W,\partial_0W;Z[\pi_1N])@>>>C_{*-1}(\partial_0W;Z[\pi_1N])\\
@VVV @VVV @VV=V\\
C_*(N;Z[\pi_1N])@>>>C_*(N,\partial_0W;Z[\pi_1N])@>>>C_{*-1}(\partial_0W;Z[\pi_1N]),\endCD$$
and similarly for $(K,\partial_0W)\to (N,\partial_0W)$. The complex $C_*(N,\partial_0W;Z[\pi_1N])$ is the
$Z[\pi_1N]$ reduction of the stratified chains, so the commutativity hypothesis on the maps to stratified
chains shows these inclusion diagrams commute with the chain maps $d_*$ and $\hat d_*$ described above. Now
we again use the fact that $\hat d_2$ preserves bases. This and the commutativity mean a 2-cell of $K$ and
the corresponding 2-handle of $W$ have the same image in $C_2(N;Z[\pi_1N])$. This in turn implies that the
difference between them represents the trivial class in $H_2(N;Z[\pi_1N])$, as required.
\subsubhead{5.3.3  Deform 2-cells}\endsubsubhead
The situation is that we have $Z[\Gamma]$ homologies (in the common 1-skeleton) between attaching maps of
2-cells of $K$ and 2-handles of $W$, and we want to $\Gamma$-deform $W$ to $W'$ so the corresponding
homologies are realized by homotopies. Recall again that $\Gamma$-deformations don't change any of the
skeleton or chain data. Also note that a single deformation move changes only one 2-handle, so it is
sufficient to improve 2-handles one at a time.

The remainder of the proof is the same as the proof in \S7 of Part I. In brief we represent the homology by
a map of a surface into  the 1-skeleton. Since we are using (at least) $Z[\pi_1]$ coefficients this surface
factors through the universal cover of $N$. This means it extends to a map of a capped surface into $N$. We 
manipulate this capped surface until we get the geometric data of
\S3 for a deformation move. The fact that the surface comes from $Z[\Gamma]$ homology means that it actually
factors through the $\Gamma$ cover, but this  plays no role in the construction until the very end when it
is used to recognize the move as a $\Gamma$ deformation. 

We warn again that in the proof of \S7 of Part I, the notations $M$ and $W$ are interchanged from their
meanings here.

Finally we have the grope refinement, extending 2.3.
\proclaim{5.4 Theorem (Relative $\Gamma_n$-Spine Realization)} Suppose $N=M\cup W$ is a $\Gamma$
decomposition and
$(K,\partial_0W)\to (N,\partial_0W)$ is a relative 2-complex as in 5.3, but with  a factorization
$\pi_1K\to
\Gamma_n(W)\to
\pi_1N$. Then there is a $\Gamma$ grope deformation of height $n+1$ from $M$ to a $\Gamma$ decomposition
$M'\cup W'$ so that the spine of $(W',\partial_0W)$ realizes $(K,\partial_0W)$ with its $\Gamma_n$
factorization, if and only if the chain hypotheses of 5.3 hold over $\Gamma_n(W)$.
\endproclaim
\demo{Proof}
Follow the proof of 5.3, replacing $\Gamma$ by $\Gamma_n(W)$, up to the point in 5.3.3 where the homology
data is used to get a surface  in the $\Gamma_n(W)$ cover of the 1-skeleton giving a homology between
attaching maps. Choose standard basis curves for the homology of this surface, but now observe that the
lift to $\Gamma_n(W)$ implies that these curves bound uncapped gropes of height $n$ mapping into the
$\Gamma$ cover of $W$. Adding these to the original surface gives an uncapped grope of height $n+1$. As
before the factorization through $\Gamma\to \pi_1N$ implies that the uncapped grope in $W$ extends to a
map of a capped grope into $N$. 

From here the proof again essentially follows Part I \S7: we see that by changing the attaching maps in
$K$ by homotopy we can reduce the singularities in the map of the capped grope until we get embedded data
as in 3.5, and can then do a grope deformation move on a handle of $W$. The surface and grope cases are
not essentially different so we will not give details. See Conant and Teichner \cite{CT} for similar
arguments for desingularizing maps of gropes. 
\enddemo

\head 6. The solvable tower\endhead
The chain invariant of Section 4  is intended to be used to distinguish different decompositions
of a 4-manifold. There is a difficulty in that  before the invariant is defined we need to know there are
some similarities, particularly in fundamental groups. The ``solvable tower'' is based on the fact that some fundamental group similarities can be detected
homologically
because the abelianization of $\pi_1$  is $H_1$. The outcome is an ``obstruction theory'': a sequence of invariants so that at each level
either the invariant distinguishes the decompositions, or the next level invariant is defined.

As before we fix a smooth compact 4-manifold $N$, a homomorphism $\gamma\:\Gamma\to\pi_1N$, and a
decomposition $\partial N=\partial_0M\cup \partial_0W$ with a factorization $\pi_1\partial_0W\to \Gamma\to
\pi_1N$. We will be considering derived groups $\Gamma_n(W)$, so nontrivial examples may arise even if we
start with $\Gamma =\pi_1N$. 

\subhead 6.1 Homological equivalence\endsubhead Suppose $M\cup W$ and $M'\cup W'$ are
$\Gamma$ decompositions of $N$. We define a {\it homological $\Gamma_n$ equivalence\/} between them to be
$(\phi,f)$, where
\roster\item $\phi$ is an isomorphism $\Gamma_n(W)\to \Gamma_n(W')$ that commutes with the
homomorphisms from $\pi_1\partial_0W$ and to
$\Gamma$; and
\item $f$ is a simple chain homotopy equivalence of $Z[\Gamma_n]$ complexes,
$$C_*(W,\partial_0W;Z[\Gamma_n])@>{f}>> C_*(W,\partial_0W;Z[\Gamma_n])$$
that chain homotopy commutes with the maps to the stratified chains $C_*(N,\partial_0W;Z[\gamma_n])$.
\endroster
In (2) we are using the isomorphism $\phi$ to identify the $\Gamma_n$ groups, and $\gamma_n\:\Gamma_n\to
\pi_1N$ is the natural map. 
\proclaim{6.2 Lemma} Suppose $(\phi,f)$ is a $\Gamma_n$ homological equivalence of decompositions. Then
$\phi$ lifts to a canonical isomorphism
$\hat\phi\:\Gamma_{n+1}(W)\to \Gamma_{n+1}(W')$ and $f$ lifts to a $Z[\Gamma_{n+1}]$  equivalence 
 $\hat f\:C_*(W,\partial_0W;Z[\gamma_{n+1,n}])\to C_*(W',\partial_0W;Z[\gamma_{n+1,n}])$ of
stratified chains.
\endproclaim
The   coefficient homomorphism in the  stratified chain groups is the quotient
$\gamma_{n+1,n}\:\Gamma_{n+1}(W)\to
\Gamma_{n}(W)$.
\demo{Proof}
It is sufficient to do this for $\Gamma_1$, since $\Gamma_{n+1}=(\Gamma_n)_1$, and the notation is easier.
In this case $\phi$ is the identity $\Gamma\to \Gamma$, and $f$ is a simple equivalence of $Z[\Gamma]$
chains.

The basic idea is that the kernel of $\Gamma_1(W)\to\Gamma$ is $H_1(W;Z[\Gamma])$. If $W$ and $W'$ have the
same $\Gamma$, and the same $H_1$ by homological equivalence, then they should have the same $\Gamma_1$.
There are two issues to deal with:  in the relative case we get a relative $H_1$ rather than the
absolute one needed, and there is a group extension problem.

By hypothesis the chain map $f$ homotopy commutes with the maps to the stratified chains
$C_*(N,\partial_0W;Z[\gamma])$. The compositions with the boundary map to $C_{*-1}(\partial_0W;Z[\Gamma])$
therefore also commute. This shows that the right-hand square in the diagram below homotopy-commutes, and
since the rows are exact this induces the chain map $g$:
$$\CD C_*(W;Z[\Gamma])@>>>C_*(W,\partial_0W;Z[\Gamma])@>>> C_{*-1}(\partial_0W;Z[\Gamma])\\
@VV{g}V @VV{f}V @VV{=}V\\
C_*(W';Z[\Gamma])@>>>C_*(W',\partial_0W;Z[\Gamma])@>>> C_{*-1}(\partial_0W;Z[\Gamma])\endCD$$

The chain map $g$ is a simple equivalence since the other two vertical maps are. According to the
1-skeleton alignment lemma we can deform the spines of $W$ and $W'$ to have the same 1-skeleton. Denote the 
fundamental group of this 1-skeleton by $F$. The kernel of $F\to\Gamma$ is the fundamental group of
1-skeleta for both
$\Gamma$ covers; denote this kernel by $K$. By abelianization $K$ maps to $H_1(W,Z[\Gamma])$ and  
$H_1(W',Z[\Gamma])$. The chain equivalence gives an isomorphism of these groups, and since the chain map is
the identity in degree 1 the maps
$K\to H_1$ commute with this isomorphism. This means their kernels are the same. But
$\Gamma_1(W)=F/\text{ker}(K\to H_1(W;Z[\Gamma])\,)$ and similarly for $\Gamma_1(W')$. Therefore the kernels
being the same gives an identification $\hat\phi\:\Gamma_1(W)= \Gamma_1(W')$.

To complete the lemma we must lift $f$ to a chain equivalence of stratified chains. The stratified chains
of $W$ are defined to be the pullback in the diagram
$$\CD C_*(W,\partial_0W;Z[\gamma_{1,0}])@>>>C_{*-1}(\partial_0W;Z[\Gamma_1(W)])\\
@VVV@VVV\\
C_*(W,\partial_0W;Z[\Gamma])@>>>C_{*-1}(\partial_0W;Z[\Gamma])\endCD$$
and similarly for $W'$.  To get a chain equivalence of stratified chains it is sufficient to get an
equivalence of pullback data for $W$ and $W'$. On the lower left we have the chain equivalence $f$. On the
lower right we have the identity, or more precisely the equivalence induced by the identity on
$\partial_0W$ and the isomorphism $\phi\:\Gamma = \Gamma$. On the upper right we have a similar
equivalence, induced by the identity on $\partial_0W$ and the new isomorphism $\hat\phi\:\Gamma_1(W)\to
\Gamma_1(W')$. These are compatible, so by universality of the pullback construction they define a chain
equivalence of stratified chains. 
\enddemo

\subhead 6.3 The chain invariant of a homological equivalence\endsubhead
Suppose $(\phi,f)$ is a homological $\Gamma_n$ equivalence from $M\cup W$ to $M'\cup W'$. Then the
{\it chain invariant\/} $\tau_{n+1}(\phi,f)$ is defined to be the equivalence class of
$(C_*(W',\partial_0W;Z[\Gamma_{n+1}]), q)$ in $\sc(Z[\Gamma_{n+1}])\downarrow
C_*(W,\partial_0W;Z[\gamma_{n+1,n}])$. 
\medskip
Here $\Gamma_{i}$ refers to either $W$ or $W'$, with the two groups being identified via $\phi$ when $i=n$,
and by the lift $\hat\phi$ of 6.2 when $i=n+1$. $\gamma_{n+1,n}$ is the natural quotient 
$\Gamma_{n+1}\to \Gamma_n$, and $q$ is the $Z[\Gamma_{n+1}]$ chain map
obtained by composing the change-of-coefficient map with the
equivalence of 6.2:
$$C_*(W',\partial_0W;Z[\Gamma_{n+1}])@>>>C_*(W',\partial_0W;Z[\gamma_{n+1,n}])@>\hat
f>>C_*(W,\partial_0W;Z[\gamma_{n+1,n}]).$$
Recall (4.3) that $\sc(Z[\Gamma_{n+1}])\downarrow D_*$ is defined when $D_*$ is a $Z[\Gamma_{n+1}]$ complex
whose $Z[\Gamma_{n}]$ reduction is a small chain object. This set is the equivalence classes of $(E_*,q)$
where $E_*$ is a small chain object over $Z[\Gamma_{n+1}]$ and $q\:E_*\to D_*$ is a chain map whose
$Z[\Gamma_{n}]$ reduction is a simple equivalence.

Note that the identity homomorphism on $\Gamma_{n}(W)$ and the identity map on chains defines a
homological equivalence of $M\cup W$ to itself. This has a chain invariant, namely
$(C_*(W,\partial_0W;Z[\Gamma_{n+1}]),q)$ where $q$ is the change-of-coefficient map to the stratified
chains. We will say the chain invariant of
$(\phi,f)$ is {\it trivial\/} if it is the same as the invariant of the identity. Note the sets in
which the chain invariants take values have not been given group structures, and even if they have group
structures ``trivial'' may not be the same as ``=\,0''. It seems reasonable to hope that these sets have
affine structures (free effective action of an abelian group) so choosing a basepoint gives an abelian group
structure with this point as 0. In this case choose the invariant of the identity as the basepoint, and
``trivial'' does become the same as ``=\,0''. The $K$-theory speculation of 5.3 is in part an effort to
find an affine structure.

\proclaim{6.4 Theorem} 
Suppose $M\cup W$ is a $\Gamma$ decomposition of $N$. 
\roster\item Every element of $\sc(Z[\Gamma_{n+1}])\downarrow
C_*(W,\partial_0W;Z[\gamma_{n+1,n}])$ is realized as $\tau_{n+1}$ of a homological $\Gamma_n$ equivalence;
\item $\tau_{n+1}$ of a homological $\Gamma_n$ equivalence is unchanged by composition with a homological
$\Gamma_{n+1}$ equivalence, so in particular is unchanged by $\Gamma_{n+1}$-deformation, or $\Gamma$ grope
deformation of height $n+2$.
\item Suppose  $(\phi,f)$ and $(\lambda,g)$ are homological $\Gamma_n$ equivalences to
decompositions $M'\cup W'$ and $\widehat M\cup \widehat W$ respectively. Then there is a homological
$\Gamma_{n+1}$ equivalence $(\alpha, h)$ from $M'\cup W'$ to $\widehat M\cup \widehat W$ with
$\alpha$ the canonical isomorphism on $\Gamma_{n+1}$ and $h$ a chain equivalence lifting the
$\Gamma_n$ equivalence $g\,f^{-1}$, if and only if 
$\tau_{n+1}(\phi,f)=\tau_{n+1}(\lambda,g)$.
\endroster
\endproclaim
\demo{Proof} This follows easily from previous results and the definitions. Statement (1) comes from 4.4.
The first part of statement (2) is a consequence of the equivalence relation used in the definition of 
$\sc(Z[\Gamma_{n+1}])\downarrow C_*(W,\partial_0W;Z[\gamma_{n+1,n}])$, and the second part comes from the
fact (3.4) that deformations induce chain equivalences. Finally (3) is again a restatement of the
definition of the equivalence relation on chain invariants.
\enddemo

\subhead 6.4 The obstruction theory\endsubhead
Suppose $M\cup W$ and $M'\cup W'$ are $\Gamma$-decompositions of $N$. Theorem 6.3 gives a plan for
comparing them: first see if they are homologically $\Gamma$ equivalent. If not they are different. If so
choose an equivalence $(\gamma_0,f_0)$ and consider the chain invariant $\tau_1(\gamma_0,f_0)$. If this is
nontrivial (different from the invariant of the identity) then the decompositions are at least somewhat
different. If the invariant is trivial then there is a homological $\Gamma_1$ equivalence. Continue in the
same way: at each stage we either see some difference or can proceed to the next stage. This plan suffers
from the usual defects, which we now discuss. The hope is that the context here is sufficiently more
explicit than the classical link setting that the defects can be better analysed.
\subsubhead 6.4.1 Indeterminate \endsubsubhead
Triviality of an invariant $\tau_{n+1}$ means a
$\Gamma_n$ equivalence lifts to a $\Gamma_{n+1}$ equivalence, but does not specify a particular lift. The
invariant $\tau_{n+2}$ at the next stage  depends on the choice of lift, so is not well-defined as a
function of the original data.  The failure to be well-defined can be related to chain self-equivalences of
the chains of $(W,\partial_0W)$. 
\subsubhead 6.4.2 Inconclusive \endsubsubhead
Nontriviality of $\tau_{n+1}$ means a particular $\Gamma_n$ equivalence does not lift to a $\Gamma_{n+1}$
equivalence. It may be that some other $\Gamma_n$ equivalence  does lift. Again this can be
analysed in terms of chain self-equivalences of the chains of $(W,\partial_0W)$.
\subsubhead 6.4.3 Limit problems \endsubsubhead
This scheme offers a way to compare $\Gamma_n$ covers of $\Gamma$ decompositions, for arbitrarily large
$n$. We would like to be able to say that vanishing of all  these invariants gives a homological
equivalence of $\Gamma_{\infty}$ covers. The dependence on choices discussed above make even the approach
unattractive, since ``vanishing of all invariants'' must be interpreted as ``there is an infinite sequence
of choices, each of which makes the next possible.'' But worse than that, the fact that there are choices
means even if such a sequence is found there may be ``lim$_1$'' problems in fitting them together
sufficiently well to pass to a limit. 

\subsubhead 6.4.4 Special cases \endsubsubhead
It may be possible to analyse chain self-equivalences, and therefore the problems above, in special
cases where activity is largely confined to a single degree. Interesting examples are suggested in~5.3.2.


\Refs
\widestnumber\key{COT}

\ref\key CM\by Tim Cochran and Paul Melvin \paper Finite type invariants of 3-manifolds  \jour Invent. Math. \vol 140 \yr 2000\pages 45--100\paperinfo
math.GT/9805026\endref

\ref\key COT\by Tim Cochran, Kent Orr, and Peter Teichner\paper Knot concordance, Whitney towers, and $L^2$
signatures\jour Ann. of Math. \vol 157 \yr 2003\pages 433--519\paperinfo math.GT/9908117\endref

\ref\key CT\by James Conant and Peter Teichner\paper Grope cobordism of classical knots\paperinfo Preprint
January 2001
\endref

\ref\key FQ \by Michael Freedman and Frank Quinn\book Topology of 4-manifolds\yr 1990\publ Princeton
University Press\endref

\ref\key Q\by Frank Quinn\paper Dual decompositions of 4-manifolds\jour Proc. Amer. Math.
Soc.\vol 354\yr 2002\pages 1373--1392\endref

\ref\key R\by Andrew Ranicki\paper Noncommutative localization in topology\toappear\endref


\endRefs
\bye